\documentclass{amsart}

\usepackage{txfonts}
\usepackage{amsmath}
\usepackage{amssymb}
\usepackage{cite}
\usepackage{graphicx, epsfig}
\usepackage{hyperref}

\usepackage[margin=3cm]{geometry}

\newtheorem{theorem}{Theorem}[section]
\theoremstyle{plain}

\newtheorem{definition}[theorem]{Definition}

\newtheorem{lemma}[theorem]{Lemma}

\newtheorem{proposition}[theorem]{Proposition}
\newtheorem{remark}[theorem]{Remark}

\numberwithin{equation}{section}

\usepackage{hyperref}
\hypersetup{
    colorlinks   = true,
    urlcolor     = blue,
    citecolor    = red,
    bookmarksopen=true
}

\title[Regularity for Elliptic Equations]{Regularity for Fully Nonlinear Elliptic Equations with Natural Growth in Gradient and Singular Nonlinearity}
\author{Mohan Mallick}
\address[Mohan Mallick]{Department of Mathematics, Visvesvaraya National Institute of Technology Nagpur, India-440010}
\email{mohan.math09@gmail.com, mohanmallick@mth.vnit.ac.in}
\author{Ram Baran Verma}
\address[Ram Baran Verma]{SRM University Amaravati, Andhra Pradesh-522502, India}
\email{rambaran.v@srmap.edu.in, rambv88@gmail.com}
\date{\today} 
\begin{document}

\maketitle

\begin{abstract}
\noindent In this article we consider the following boundary value problem
 \begin{equation*}\label{abs}
\left\{
\begin{aligned}
F(x,u,Du,D^{2}u)+c(x)u+ p(x)u^{-\alpha}&=0~\text{in}~\Omega\\
u&=0~~\text{on}~~\partial\Omega,
\end{aligned}
\right.
\end{equation*}
where $\Omega$ is a bounded and $C^{2}$ smooth domain in $\mathbb{R}^N$ and $F$ has superlinear growth in gradient and $c(c)<-c_{0}$ for some positive constant $c_{0}.$ Here, we studies the boundary behaviour of the solutions to above equation and establishes the global regularity result similar to one established in \cite{felmer2012existence,gui1993regularity} with linear growth in gradient.  
\end{abstract}
\section{Introduction}
Consider the following boundary value problem:
\begin{equation}\label{main}
\left\{
\begin{aligned}
F(x,u,Du,D^{2}u)+c(x)u&=-p(x)u^{-\alpha}~\text{in}~\Omega\\
u&=0~~\text{on}~~\partial\Omega,
\end{aligned}
\right.
\end{equation}
where $\Omega$ is a $C^{2}$ smooth domain in $R^{N},$ $F$ has superlinear growth in gradient and $p$ is a positive function. There are many problems in sciences and engineering the study of which leads to an equation of the form \eqref{main}. We would like mentioned \cite{bertozzi1998long,bertozzi2000finite,pelesko2002mathematical}, where first two articles deal with the study of steady state of thin
films while third is related to the modelling of the MEMS devices. The study of existence, boundary behaviour and regularity of solutions to \eqref{main} in the case where $F$ is a linear elliptic operator has linear growth in gradient has seek attention of a large group of researcher. The first attempt to establish the existence and understand the behaviour of solution to equations having singular non linearity has been made in \cite{crandall1977dirichlet}. In this article, one of the results asserts that in general these problems do not have classical solutions. Later on A. C. Lazer and P. J. Mckenna in \cite{lazer1991singular},consider \eqref{main} with $F=\Delta$ and $c=0.$ Here authors proved that if the domain is smooth and $p$ is H\"{o}lder continuous function then there exists a classical solution and is comparable to the power of the eigen function of the associated Dirichlet problem. In this continuation, Gomes \cite{gomes1986singular} has considered the \eqref{main} the with $F=\text{Div}\Big(A\nabla u\Big),$ $c=0$ and $p$ behaving like a power of the distance function and studied the behaviour and regularity of solution. In 1993, C. Gui and Lin \cite{gui1993regularity} consider the \eqref{main} with Laplace operator, $c=0$ and $p$ behaving like \eqref{est1} and prove that the solution is H\"{o}lder continuous or gradient of the solution is H\"{o}lder continuous depending on the range of $\mu.$ Recently, in \cite{felmer2012existence} P. Felmer etc. consider \eqref{main} with $F$ having linear growth in gradient and prove the existence of positive solution and studies its regularity properties. The existence of solution to \eqref{main}  with sublinear nonlinear term as well as  associated semipositone problems have been established in \cite{tyagi2017positive}. Although the existence of solution to a class of equations having superlinear growth in gradient and singular non linear term is proved in \cite{tyagi2019positive} but the regularity of the solutions to such equations has not been proved  yet.\\
The equations having superlinear growth in gradient have been attaracting continuous attention of the researchers. Study of such equations was started by L. Boccardo et.al\cite{boccardo1982existence} in 1982. Since then many authors considered the problems involving natural growth in gradient, for example, see\cite{boccardo1983existence,bensoussan1988non,maderna1992quasilinear,ferone2000nonlinear, boccardo1992estimate}. Fully nonlinear elliptic equation with natural growth in gradient first of all appears in \cite{ishii1990viscosity,crandall1992user}. In \cite{ishii1990viscosity},  authors establish the comparison principle and existence of solutions. A comprehensive study of this class of equations have been by B. Sirakov in \cite{sirakovsolvability}. One of the results of this article asserts global H\"{o}lder estimate for solutions see Theorem 2 \cite{sirakov2010}. At this point we would like to point out that we do not have Harnack inequality for this class of equation in general. B. Sirakov uses a version of weak(very weak see Remark below Proposition 4.2\cite{sirakov2010}). Recently, G. Nornberg \cite{nornberg2019c1} establishes $C^{1,\alpha}$ regularity for these class of equations \cite{nornberg2019c1}. This result asserts the interior regularity of solutions \eqref{main}. Here we establish the global regularity of solutions \eqref{main}. In this article we establish the regularity of solutions to \eqref{main} similar to one established in \cite{felmer2012existence}. As in \cite{felmer2012existence} we first prove the boundary behaviour of solutions and then establish the regularity of solution. These boundary behaviour and regularity of solutions depends on the range of $\mu$ as discussed above. For a range of values of $\mu$ we prove global H\"{o}lder continuity of the gradient. This result has been proved by following the Krylov approach. For linear equation without lower order we refer Theorem 9.31\cite{GT} while for nonlinear equation linear growth in gradient we refer \cite{silvestre2014boundary}. We also refer to \cite{adimurthi2020twice} for parabolic analogue of this result. Recently, this result has been generalised to the equations with superlinear growth in gradient \cite{braga2020krylov}. But the approach used in this article requires the non-homogeneous term to be in $L^{N}$ which may not be the case in our case. We follow the approach used in \cite{GT} with proper adaptation to the equations with singular non linearity introduced in \cite{felmer2012existence}. The problem considered here has superlinear growth in gradient in gradient as well as singular nonlinearity. In \cite{felmer2012existence} authors considered equation with linear growth in gradient and Harnack inequality is available for these class of equations. However, Harnack inequality is not available for the class of equations considered here. Our proof uses only weak Harnack inequality for equation with linear growth in gradient. At this point we would like bring the notice that weak Harnack inequality true for any nonnegative supersolution. However, for Harnack inequality to be true the function has to be non negative sub and super solution of the extremal equations. This is very important for our result because we apply the transformation \eqref{trans} to get rid of the superlinear growth in gradient but we get only on side inequality. This approach requires the global Lipschitz continuity of the solution. This is the content of our next section. Moreover, in section two we prove the boundary behaviour for the positive solutions to equations of \eqref{main}. All the result of this section requires the comparison principle to hold for the class of equations considered. Our operator $F$ is proper as well as $c(x)<-c_{0}$ guarantees that the comparison principle hold for the details see, Proposition 3.1\cite{sirakovsolvability}.
\section{Definition and Basic Results}
The function $F$ in \eqref{main} has superlinear growth in gradient and continuous in its argument $x.$ Moreover, it satisfies the following structure condition:
\begin{equation}
(SC)\left\{
\begin{aligned}
&(i)~~\mathcal{M}_{\lambda, \Lambda}^{-}(M_{1}-M_{2})-B(|p_{1}|+|p_{2}|)|p_{1}-p_{2}|-b|p_{1}-p_{2}|-d((r_{1}-r_{2})^{+}) \\
&\leqq F(M_{1},p_{1},r_{1},x)-F(M_{2},p_{2},r_{2},x)\\
&\leq\mathcal{M}_{\lambda, \Lambda}^{+}(M_{1}-M_{2})+B(|p_{1}|+|p_{2}|)|p_{1}-p_{2}|+b|p_{1}-p_{2}|+d((r_{2}-r_{1})^{+}),\\
&(ii)~~F(x,0,0,0)=0.
\end{aligned}
\right.
\end{equation}
where $B,b$ and $d$ are positive constants and $(x,r_{1},p_{1}, M_{1}),~(x,r_{2},p_{2},M_{2})\in\Omega\times\mathbb{R}\times \mathbb{R}^{N}\times S(N),$ where $S(N)$ denotes the set of real symmetric matrices of order $N\times N.$ Associated to the constant $B$ appearing in the above $(SC)$ we define the following two constants:
\begin{equation}\label{m}
l_{1}=B/\Lambda~~~\text{and}~~l_{2}=B/\lambda.
\end{equation}
Observe that $l_{1}<l_{2}$ as long as $\lambda<\Lambda.$ The operator $F$ is proper uniformly elliptic operator. As we mentioned that this work is a generalization of the works \cite{sirakovsolvability,nornberg2019c1,felmer2012existence,braga2020krylov} so we will be borrowing the required results from these articles. Since many estimates near the boundary have been obtained by construction the barriers, so we have taken the coefficients continuous and bounded. Consequently by solution in this article we wean continuous "viscosity solution" unless otherwise explicitly mentioned. Consider the following problem
\begin{equation}\label{ref}
\left\{
\begin{aligned}
F(x,u,Du,D^{2}u)+c(x)u&=g(x,u)~\text{in}~\Omega\\
u&=0~~\text{on}~~\partial\Omega,
\end{aligned}
\right.
\end{equation}
where $g~:~\Omega\times\mathbb{R}\longrightarrow\mathbb{R}$ is a continuous function.
\begin{definition}[see \cite{MHP}]\label{viscosity}
A function $u\in C(\bar{\Omega})$ is called the viscosity subsolution (resp., supersolution) of \eqref{ref} if for any $\phi\in C^{2}(\Omega)$ such that $u-\phi$ has a maximum (resp., minimum) at some point $x_{0}\in \Omega$ then we have
\begin{equation*}
F(x_{0},u(x_{0}), D\phi(x_{0}),D^{2}\phi(x_{0}))+c(x)u\geq g(x_{0},u(x_{0}))~~(\text{resp.,}~~\leq g(x_{0},u(x_{0}))).
\end{equation*}
$u$ is a solution if it is both subsolution and supersolution at the same time and $u=0$ on $\partial\Omega.$
\end{definition}

We will be using the following transformation again and again, for more details see Lemma 2.3 \cite{sirakovsolvability}.
\begin{lemma}\label{trans}
Let $u\in C^{2}(\Omega)$. For any $l>0$ set
\[v=\frac{e^{lu}-1}{l},~~\text{and}~~w=\frac{1-e^{-lu}}{l}\]
Then in $\Omega$ we have \[Dv=(1+lv)Du,~~\text{and}~~Dw=(1-lw)Du\] and 
\[l\lambda|Du|^2+\mathcal{M}_{\lambda,\Lambda}^{\pm}(D^2u)\leq\frac{\mathcal{M}_{\lambda,\Lambda}^{\pm}(D^2v)}{1+lv}\leq l\Lambda|Du|^2+\mathcal{M}_{\lambda,\Lambda}^{\pm}(D^2u)\]
\[-l\Lambda |Du|^{2}+\mathcal{M}_{\lambda,\Lambda}^{\pm}(D^2u)\leq\frac{\mathcal{M}_{\lambda,\Lambda}^{\pm}(D^2w)}{1-lw}\leq-l\lambda |Du|^{2}+\mathcal{M}_{\lambda,\Lambda}^{\pm}(D^2u) \]
and, clearly, $u=0$(resp. $u>0$) is equivalent to $v=0$(resp. $v>0$). This lemma is also remains true if $u\in C(\Omega)$ and respective classical inequalities hold in the viscosity sense.
\end{lemma}
Above lemma helps us to absorb the superlinear growth in the gradient but at the same time dependent variable gets multiplied by a term containing  the logarithmic of the function. In order to deal with such non-linearities we will be using the following elementary result from the calculus. 
\begin{lemma}\label{calc}
Let us consider 
\[f(r)=-(1/l_{2})\log(1-l_{2}r)-t~~~\text{and}~~g(r)=(1/l_{1})\log(1+l_{1}r)-r\]
defined respectively in $[0,1/l_{2})$ and $[0,\infty).$ The function $f$ is non-decreasing and has minimum at $r=0.$ While the function $g$ is non-increasing and has maximum at $r=0.$ In particular, we have 
\begin{equation*}\label{enarm}
\begin{aligned}
&(i)~-(1/l_{2})\log(1-l_{2}r)\geq r~~\text{for}~~t\in[0,1/l_{2})]\\
&(ii)~(1/l_{1})\log(1+l_{1}r)\leq r~~\text{for}~~t\in[0,\infty).
\end{aligned}
\end{equation*}
\end{lemma}
As we have mentioned in the introduction that this article discusses about the global regularity of solutions of \eqref{main} depending on the range of $\alpha.$ One of our results is about the global $C^{1,\beta}$ regularity of solutions. This regularity has been established as usual by combining the interior gradient estimate in \cite{nornberg2019c1} and boundary H\"{o}lder continuity in the standard way. The boundary gradient H\"{o}lder continuity has been obtained by following Krylov approach. For more details about this approach for obtaining the gradient estimate at boundary, we refer to Theorem 9.31 \cite{GT} and \cite{braga2020krylov,silvestre2014boundary}. At this point, we would like to point out that even though this result has been established in \cite{braga2020krylov} for more general class of operators but it has been assumed that the non-homogeneous term is in $L^{N}$ which may not be true in our case. We have obtained this result by following the approach as in \cite{felmer2012existence} which uses modified auxiliary function suiting to equation with singular nonlinear term. But this approach require Harnack inequality which is not available for the considered class of equations. We have established this result by using the weak Harnack inequality followed by transformation in Lemma\ref{trans}. For more details we refer to the proof of Proposition \ref{propkrylov} and  Theorem \ref{bmain}. We have taken the following weak Harnack inequality from \cite{N} see Proposition 1.8 there. 
\begin{theorem}\label{weakh}
Suppose that $u$ is a non negative viscosity supersolution of 
\begin{equation*}
\mathcal{M}_{\lambda, \Lambda}^{-}(D^{2}u)-b|Du|\leq f~~\text{in}~~B_{2R}(x_{0}).
\end{equation*}
Then then there exist $p$ and $C$ depending on $\lambda, \Lambda,n, b$ such that 
\begin{equation*}\label{enarm}
\begin{aligned}
&\Bigg[\fint_{B_{R}}u^{p}dx\Bigg]^{1/p}\leq C\Bigg(\inf_{\substack{B_{R}}}u+R\|f\|_{L^{N}(B_{R})}\Bigg)
\end{aligned}
\end{equation*}
holds.
\end{theorem}
As we know that the krylov approach requires the boundary Lipschitz estimate for the solution to hold. This estimate as in \cite{felmer2012existence}, will be obtained as by constructing an appropriate sub and supersolution. These sub and supersolution are solution to a boundary value problem. We establish the existence of solution to the above mentioned auxiliary class of equations by monotone iteration method. As we know that this method require the construction of a subsolution and a supersolution along with the comparison principle. In the construction of sub and supersolution as in Lemma 1\cite{felmer2012existence}, we need the eigen function of certian equation. We start by introduction of the following two operators:
\begin{equation}
F_{1}^{\pm}(u)=\mathcal{M}_{\lambda,\Lambda}^{\pm}(D^2u)\pm B|Du|^{2}\pm b|Du|
\end{equation}
\begin{equation}
F_{2}^{\pm}(u)=\mathcal{M}_{\lambda,\Lambda}^{\pm}(D^2u)+\pm b|Du|.
\end{equation}
For given $0<r<R$ we define the annular domain
\begin{equation}
C_{r, R}=\{x\in\mathbb{R}^{N}~|~~r<|x|<R\}.
\end{equation}
For the following existence result see Theorem 1.1 \cite{AB}.
\begin{theorem}
There exists $(\lambda^{+}_{1},\psi)\in\mathbb{R}^{+}\times\Big(W^{2,N}_{loc}(C_{r,R})\cap C(\overline{ {C_{r,R}}} )\Big)$ solution to the following 
\begin{equation}\label{B4}
\left\{
\begin{aligned}
F_{2}^{+}\left(\psi\right)&=-\lambda^{+}_{1}\psi~~\text{in} ~~C_{r, R},\\
\psi&>0~~\text{in}~~C_{r, R}\\
\psi&=0~~\text{on}~~\partial C_{r, R}.\\
\end{aligned}
\right.
\end{equation}
\end{theorem}
We also need the following result which asserts the existence of positive solution to the weighted eigenvalue problem for details see Theorem 6\cite{felmer2012existence}.
\begin{theorem}
There exists $(\lambda^{+},\phi)\in\mathbb{R}^{+}\times W^{2,N}_{loc}(A_{r,R})$ such that 
\begin{equation}\label{B4w}
\left\{
\begin{aligned}
F_{2}^{+}\left(\phi\right)&=-\lambda^{+}(\gamma-r)^{\mu}\phi~~\text{in} ~~C_{r, R},\\
\phi&>0~~\text{in}~~C_{r, R}\\
\phi&=0~~\text{on}~~\partial C_{r, R}.\\
\end{aligned}
\right.
\end{equation}
\end{theorem}
We also need the following comparison principle for more details see, Proposition 3.1\cite{sirakovsolvability}.
\begin{proposition}\label{comp1}
Let $F$ satisfies $(SC)$ and $u$ and $v$ be respectively sub and supersolution of
\begin{equation}\label{B4c}
\left\{
\begin{aligned}
F(x,u, Du, D^{2}u)-c_{0}u&\geq f(x)~~\text{in} ~~\Omega,\\
F(x,v, Dv, D^{2}v)-c_{0}v&\leq f(x)~~\text{in} ~~\Omega,\\
\end{aligned}
\right.
\end{equation}
for some $c_{0}>0$ and $u\leq v$ on $\partial\Omega$ then $u\leq v$ in $\Omega.$
\end{proposition}
\section{Estimates on the solution near the boundary of the domain}
 Along with this establish the boundary behaviour of solutions to our class of equation will also be accomplished by establishing the existence and properties of the solutions to an appropriate class of equations.


\begin{lemma}
Assume $\mu \geq 0, \alpha>0$. There exists $R_0>0$ such that for every $0<r<$ $R \leq R_0, C>0$ and $M>0$, the problem
\begin{equation}\label{ra1}
\begin{aligned}
F^{1}_{+}(v)-c_{0}u+M(r-\varrho)^\mu v^{-\alpha}&=0 \text {in} C_{\varrho, R}, \\
v&=0 \text { on } \partial B_\varrho,\\
\quad v&=L~\text{on}~\partial B_R,
\end{aligned}
\end{equation}
has unique positive solution.
\end{lemma}
\begin{proof}
Let us consider
\begin{equation}\label{B66}
F^{1}_{+}(w)-c_{0}w+\frac{M(r-\varrho)^\mu}{(w+\delta)^{\alpha}}=0
\end{equation}
regularized version of \eqref{ra1}. In view of the comparison principle Proposition \ref{comp1}, it is sufficient to construct an appropriate ordered sub and supersolution. Let us proceed to construct the subsolution. Fix a $1<\kappa<2,$ and consider $u_{1}=m_{1}\psi^{\kappa}.$ Then 
\begin{equation*}\label{B65}
\begin{aligned}
&F^{2}_{+}(u_{1})-c_{0}u_{1}+\frac{M(r-\varrho)^\mu}{(u_{1}+\delta)^{\alpha}}\geq m_{1} \psi^{\kappa}\Big[\kappa(\kappa-1)\lambda|\psi^{\prime}|^{2}\psi^{-2}-c_{0}\Big]+(r-\varrho)^{\mu}\Big[-m_{1}\kappa\lambda^{+}\psi^{\kappa}+\frac{M}{(m_{1}\psi^{\kappa}+\delta)^{\alpha}}\Big].
\end{aligned}
\end{equation*}
Note that in a neighbourhood of $\partial C_{\varrho,R}$ first term in the right is non negative for any $m_{1}.$ In the rest part of $C_{\varrho,R}$ we can choose $m_{1}$ and $\delta_{0}$ sufficiently small such that the right hand side is nonnegative for any $\delta\in(0,\delta_{0}).$ Therefore $u_{1}=m_{1}\psi^{\kappa}$ is a subsolution of \eqref{B66} for any $\delta\in(0,\delta_{0}).$ Furthermore, it will stay a subsolution for any smaller $m_{1}.$ 
In order to construct a supersolution we first construct a supersolution of the following problem
\begin{equation}\label{B5111}
\begin{aligned}
F_{2}^{+}\left(w\right)-\frac{c_{0}}{l_{1}}(1+l_{1}w)\log(1+l_{1}w)+M(r-\varrho)^{\mu}\Big[1+l_{1}w\Big]\Big[\frac{1}{l_{1}}\log(1+l_{1}w)\Big]^{-\alpha}&\leq0~~\text{in} ~~C_{\varrho,R}.\\
\end{aligned}
\end{equation}
Let us consider the function 
\begin{equation}
u_{2}=(1+m_{2})\psi^{\eta}+m_{2}\psi+\frac{L(r-\varrho)}{(R-\varrho)}
\end{equation}
where $\psi$ is from \eqref{B4} and $m_{2}$ is a constant to be chosen below. Observe that $u_{2}$ is nonnegative therefore the term containing $c_{0}$ is negative for this function.
\begin{equation}\label{B6111}
(3.5)\left\{
\begin{aligned}
&F^{+}_{2}(u_{2})+M(r-\varrho)^{\mu}(u_{2}+\delta)^{-\alpha}\leq (1+m_{2})\eta(\eta-1)\lambda|\psi^{\prime}|^{2}\psi^{\eta-2}\\
&-(1+m_{2})\eta\lambda^{+}_{1}\psi^{\eta}-m_{2}\lambda^{+}_{1}\psi\nonumber+M(r-\varrho)^{\mu}\\
&\Bigg[1+l_{1}(1+m_{2})\psi^{\eta}+m_{2}l_{1}\psi+\frac{Ll_{1}(r-\varrho)}{R-\varrho}\Bigg]\Bigg[\frac{1}{l_{1}}\log\Big(1+l_{1}(1+m_{2})\psi^{\eta}+l_{1}m_{2}\psi+\frac{Ll_{1}(r-\varrho)}{(R-\varrho)}\Big)\Bigg]^{-\alpha}\\
\end{aligned}
\right.
\end{equation}
Let
\[ M^{+}_{\lambda,\Lambda}(D^{2}\mu(x))+b|D\mu(x)|\leq C_{1}~~\text{for~all}~x\in C_{\varrho,R}\]
where $\mu(x)=\frac{L(r-\varrho)}{(R-\varrho)}.$ We first choose a small neighbourhood of $\partial C_{\varrho,R}$ such that
\begin{equation*}
\left\{
\begin{aligned}
C_{1}+M(|x|-\varrho)^{\mu}\Big[1+l_{1}\psi^{\eta}+\frac{Ml_{1}(|x|-\varrho)}{(R-\varrho)}\Big]\Big[\frac{1}{l_{1}}\log\Big(1+l_{1}\psi^{\eta}+\frac{Ml_{1}(r-\varrho)}{(R-\varrho)}\Big)\Big]^{-\alpha}&<\frac{1}{2}\eta(1-\eta)\lambda|\psi^{\prime}|^{2}\psi^{\eta-2}\\
M(r-\varrho)^{\mu}\Big[2l_{1}\max\{\|\psi^{\eta}\|_{L^{\infty}},\|\psi\|_{\infty}\}\Big]\Big[\frac{1}{l_{1}}\log\Big(1+l_{1}\psi^{\eta}+\frac{Ml_{1}(r-\varrho)}{(\Gamma-\gamma)}\Big)\Big]^{-\alpha}&<\frac{1}{2}\eta(1-\eta)\lambda|\psi^{\prime}|^{2}\psi^{\eta-2}\\
\end{aligned}
\right.
\end{equation*}
Then outside of the chosen neighbourhood $\psi>\varpi$ for some positive constant $\varpi.$ Therefore we can choose $m_{2}$ large enough such that the right hand side of \eqref{B6111} is non positive. Therefore, in view of Lemma \ref{trans} we find that 
\begin{equation*}
\tilde{u}_{2}=\frac{1}{l_{1}}\log\Bigg(1+l_{1}(1+m_{2})\psi^{\eta}+m_{2}l_{1}\psi+\frac{Ll_{1}(r-\varrho)}{(R-\varrho)}\Bigg)
\end{equation*}
is a supersolution. Furthermore we can choose $m_{1}$ sufficiently small and $m_{2}$ large enough such that $u_{1}\leq \tilde{u}_{2}.$ The rest of the proof follows in the same manner as in Lemma 1\cite{felmer2012existence}.
\end{proof}
\section{Estimates on the solution near the boundary of the domain}
\begin{theorem}\label{boundary1}
Suppose that $u$ is a solution of \eqref{main}. Suppose also that there are two positive constants $C_{1},C_{2}$ such that 
\begin{equation}\label{est1}
C_{1}d^{\mu}(x)\leq p(x)\leq C_{2}d^{\mu}(x),
\end{equation}
where $d(x)$ is the distance of $x$ from the boundary. 
\begin{enumerate}
\item{}~If $\alpha<1+\mu$ then there exist constants $A_1, A_2>0$ such that

\[A_1 d(x)\leq u(x)\leq A_2 d(x), \quad x \in \Omega .\]

\item{} If $\alpha=1+\mu$ then there exist constants $a_1, a_2, D>0$ such that

\[a_1 d(x)(D-\log d(x))^{1 /(1+\alpha)} \leq u(x) \leq a_2 d(x)(D-\log d(x))^{1 /(1+\alpha)}.\]
\item{} If $\alpha>1+\mu$ then there exist constants $a_1, a_2>0$ such that
    
\[a_1 d(x)^{(\mu+2) /(1+\alpha)} \leq u(x) \leq a_2 d(x)^{(\mu+2) /(1+\alpha)}, \quad x \in \Omega .\]
\end{enumerate}
\end{theorem}
\begin{proof}
Since $\Omega$ is $C^{2}$ we can find a $\varrho>0$ such that if 
\[\Omega_{\varrho}=\{x\in \Omega~|~d(x)<\varrho\},\]
then for each $x\in\Omega_{\gamma}$ there exists $y=y(x)\in\partial\Omega$ and $z(x)\not\in\Omega$ along the normal at $y(x)$ such that $|y-z|=\varrho.$ 
Set $R=2\varrho$ and decrease $\varrho.$ Note that $B_{\varrho}(z)$ is an exterior tangent to $\partial\Omega.$ So for $x\in \Omega\cap B_{\varrho}(z)$~\eqref{est1} can be written as follow:
\[p(x)\leq C_{2}d^{\mu}(x)\leq C_{2}(|x-z|-\varrho)^{\mu}.\]
Now we choose the solution of \eqref{ra1} with $M\geq C_{2}$ and $L=\sup_{\Omega}u$ with a suitable $c_{0}$ so that comparison principle holds. Therefore by comparison principle
\[u(x)\leq w(x)=v(|x-z|).\]
The required upper bound will follow by the arbitrariness of $x$ if we show that 
\begin{equation}\label{ra2}
w(r)\leq A_{2}(r-\varrho),~~\text{for~all}~r\in[\varrho,~R]
\end{equation}
under the assumption $\mu-\alpha>-1$ for some positive constant $A_{2}.$ The inequality will be achieved by rewriting the \eqref{ra1} in more detail way. Let us introduce the function 
\[
    \theta(s)= 
\begin{cases}
    \Lambda,& \text{if}~~s\geq 0\\
    \lambda,& \text{if}~~s<0.
\end{cases}
\]
Then, since $w$ is radially symmetric, the eigenvalues of $D^2 w$ are $w^{\prime \prime}(r)$ and $w^{\prime}(r) / r$, so by definition of $M^{+}_{\lambda,\Lambda}$ we have
$$
M^{+}_{\lambda,\Lambda}\left(D^2 w\right)(r)=\theta\left(w^{\prime \prime}(r)\right) w^{\prime \prime}(r)+\theta\left(w^{\prime}(r)\right)(N-1) \frac{w^{\prime}(r)}{r} .
$$
Thus $w$ satisfies the equation
\begin{equation}\label{sabu}
\theta\left(w^{\prime \prime}(r)\right) w^{\prime \prime}+\theta\left(w^{\prime}(r)\right)(N-1) \frac{w^{\prime}}{r}+B|w^{\prime}|^{2}+b|w^{\prime}|-c_{0}w+C(r-\varrho)^\mu w^{-\alpha}=0.
\end{equation}
In order to write \eqref{sabu} let us define
$$
\begin{aligned}
\chi(r)& =\frac{\theta\left(w^{\prime}(r)\right)(N-1)}{\theta\left(w^{\prime \prime}(r)\right) r}, \\
\xi(r) & =\exp\left(\int_1^r \chi(s) d s\right) \text { and } \tilde{\xi}(r)=\frac{\xi(r)}{\theta\left(w^{\prime \prime}(r)\right)}.
\end{aligned}
$$
It is easy to observe the following inequalities
\begin{equation*}
\left\{
\begin{aligned}
& N_{+}-1 \leq \chi(r) r \leq N_{-}-1, \\
& r^{N_{-}-1} \leq \xi(r) \leq r^{N_{+}-1} \text { and } ~\frac{\xi(r)}{\Lambda} \leq \tilde{\xi}(r) \leq \frac{\xi(r)}{\lambda},
\end{aligned}
\right.
\end{equation*}
where $N_{+}=\frac{\lambda}{\Lambda}(N-1)+1$, and $N_{-}=\frac{\Lambda}{\lambda}(N-1)+1.$
In view of the above notation \eqref{sabu} can be rewritten as follows:
\begin{equation}
\left(\xi w^{\prime}\right)^{\prime}+\tilde{\xi}\big[B|w^{\prime}|^{2}+b|w^{\prime}|-c_{0}|w|+C(r-\varrho)^\mu w^{-\alpha}\big]=0.
\end{equation}

Let $\varrho_0=\sup \left\{r \in[\varrho, R] \mid w^{\prime}(s)>0, s \leq r\right\}$. In view of H\"{o}pf's lemma $\varrho_0>\varrho$ and $w^{\prime}(\varrho_{0})=0.$ Now, by observing that  $w(s)\geq a_{1}(s-\gamma).$ Integrating for $r \in\left(\varrho, \varrho_0\right)$, we find
\begin{equation*}
w^{\prime}(r)  \leq(\xi(r))^{-1} \int_r^{r_0} \tilde{\xi}(s)\left[ B~|w^{\prime}(s)|^{2}+b~|w^{\prime}(s)|+c_{0}|w| +C~a_1^{-\alpha}(s-\gamma)^{\mu-\alpha}\right] ds.
\end{equation*}
The right hand side is bounded in view of $\mu-\alpha>-1,0<c(\gamma) \leq \xi(r), \tilde{\xi}(r) \leq C$, and $w\in W^{2,N}_{loc}.$ Thus $w^{\prime}(r)$ is bounded for $r \in\left[\gamma, r_0\right]$, so the \eqref{ra2} follows.\\


\textbf{Estimate~}$\alpha=1+\mu$:\\
The proof of this part follow on the line of \cite{gui1993regularity,felmer2012existence}. As in the proof of the above part by the regularity of the domain we can find $\gamma$ such that for each $x\in\Omega_{\varrho}$ we can find $y(x)\in\partial\Omega$ and $z\not\in\Omega$ such that $|y-z|=\varrho.$ As we wish to show that
\begin{equation}\label{ca2}
u(x) \leq A_2 d(x)(D-\log d(x))^{1 /(1+\alpha)}
\end{equation}
for some constant $A_{2}.$ In view of the inequality 
\[\frac{1}{l_{1}}\log(1+\l_{1}t)\leq t~~\text{for all}~~t\geq 0,\]
it is sufficient to show that 
\begin{equation}\label{ca3}
u(r)\leq \frac{1}{l_{2}}\log(1+l_{2}u_{4}(r)),
\end{equation}
 where $u_{4}(r)=\bar{C}(r-\varrho)(D-\log(|x|-\varrho))^{1/(1+\alpha)},$ for sufficiently large $\bar{C}$ and  $D>1+\log (2 \operatorname{diam}(\Omega)).$ We will achieve this by comparison principle. Since  $(1+l_{1}u_4)(1/l_{1})\log(1+l_{1}u_{4})$ is non negative therefore by Lemma \ref{trans} if $u_{4}$ satisfies 
\begin{equation}\label{arbi1}
F^{+}_{2}(u_{4})+c_{2}(r-\varrho)^{\alpha-1}\Big[(1+l_{1}u_{4}(r))\Big]\Big[\frac{1}{l_1}\log\Big(1+l_{1}u_{4}(r)\Big)\Big]^{-\alpha}\leq 0~~\text{for~all}~~r\in[\varrho,R]
\end{equation}
for $\bar{C}$ sufficiently large  satisfying $\displaystyle{\sup_{\Omega}u(x)}\leq 1/l_{2}\log(1+l_{2}u_{4}(R)).$ Then $(1/l_{1})\log(1+l_{1}u_{4})$ becomes a supersolution of 
\begin{equation*}
F^{+}_{2}(u_{4})+B|u_{4}^{\prime}|^{2}-c_{0}u_{4}+c_{2}(r-\varrho)^{\mu}u^{-\alpha}\leq 0,
\end{equation*}
consequently the result follow from comparison principle. Now, we establish \eqref{arbi1}. Let $\varrho$ as above and consider $R=(1+\theta)\varrho$, with $\theta>0$ to be fixed later. We further assume $R=(1+\theta) \varrho<1.$ 

\
First and second derivative of $u_{4}(r)$ are given as follows:
\begin{equation*}
 \left\{
\begin{aligned}
&u_4^{\prime}(r)=\bar{C} h(r)\left(D-\log (r-\varrho)-\frac{1}{1+\alpha}\right),\\
&u_4^{\prime \prime}(r)=\frac{-\bar{C}}{(1+\alpha)(r-\varrho)} h(r)\left(1+\frac{\alpha}{1+\alpha}(D-\log (r-\varrho))^{-1}\right), \\
&\frac{c_2(r-\varrho)^\mu}{u_4^\alpha}=\frac{c_2 h(r)}{\bar{C}^\alpha(r-\varrho)},
\end{aligned}
\right.
\end{equation*}

\begin{equation*}
 \left\{
\begin{aligned}
&u_4^{\prime}(r)=\bar{C} h(r)(D-\log (r-\varrho)-1 /(1+\alpha)),\\
&u_4^{\prime \prime}(r)=\frac{-\bar{C}}{(1+\alpha)(r-\varrho)} h(r)\left(1+\frac{\alpha}{1+\alpha}(D-\log (r-\varrho))^{-1}\right), \\
&\frac{c_2(r-\varrho)^\mu}{u_4^\alpha}=\frac{c_2 h(r)}{\bar{C}^\alpha(r-\varrho)},
\end{aligned}
\right.
\end{equation*}

where we have used $\mu+1=\alpha$ and 
\begin{equation}
h(r)=(D-\log (r-\gamma))^{-\alpha /(1+\alpha)}.
\end{equation}
Some observation about the above derivative are as follows:
\begin{enumerate}
\item{}$u_4^{\prime}(r)>0,~~u_{4}^{\prime \prime}(r)<0.$
\item{} For all $r\in (\varrho, R)$ we have 
\begin{equation}
\left\{
\begin{aligned}
&\frac{\Lambda(N-1)(r-\varrho)}{r}\Big(D-\frac{1}{1+\alpha}-\log(r-\varrho)\Big)-\frac{\lambda}{1+\alpha}\Big(1+\frac{\alpha}{1+\alpha}[D-\log(r-\varrho)]^{-1}\Big)\\
&\leq-\frac{\lambda}{1+\alpha}+\Lambda(N-1)\theta\Big(D-\frac{1}{1+\alpha}-\log\theta\Big)\\
&\leq-\frac{\lambda}{2(1+\alpha)},
\end{aligned}
\right.
\end{equation} 
\end{enumerate}
provided $\theta$ is small enough. With the above chosen value of $\theta$ we have
\begin{equation}\label{B6}
\left\{
\begin{aligned}
&F^{+}_{2}(u_{4})+\frac{c_{2}(r-\varrho)^{\alpha-1}\Big[(1+l_{1}u_{4}(r))\Big]}{\Big[\frac{1}{l_1}\log\Big(1+l_{1}u_{4}(r)\Big)\Big]^{\alpha}}\leq\\
&\frac{\bar{C}h(r)}{(r-\varrho)}\Bigg[\frac{\Lambda(N-1)(r-\varrho)}{r}\Big(D-\frac{1}{1+\alpha}-\log(r-\varrho)\Big)-\frac{\lambda}{1+\alpha}\Big(1+\frac{\alpha}{1+\alpha}[D-\log(r-\varrho)]^{-1}\Big)+\\
&b(r-\varrho)\Big[D-\log(r-\varrho)-\frac{1}{1+\alpha}\Big]+\frac{c_{2}(r-\varrho)^{\alpha}\Big(1+l_{1}\bar{C}(r-\varrho)\big(D-\log(r-\varrho)\big)^{1/(1+\alpha)}\Big)}{h(r)\bar{C}\Big[\frac{1}{l_{1}}\log\Big(1+l_{1}\bar{C}(r-\varrho)\big(D-\log(r-\varrho)\big)^{1/(1+\alpha)}\Big)\Big]^{\alpha}}\Bigg]\\
&\leq\frac{\bar{C}h(r)}{(r-\varrho)}\Bigg[-\frac{\lambda}{2(1+\alpha)}+b(r-\varrho)\Big[D-\log(r-\varrho)-\frac{1}{1+\alpha}\Big]+\frac{c_{2}(r-\varrho)^{\alpha}\Big(1+l_{1}(r-\varrho)\bar{C}\big(D-\log(r-\varrho)\big)^{1/(1+\alpha)}\Big)}{h(r)\bar{C}\Big[\frac{1}{l_{1}}\log\Big(1+l_{1}\bar{C}(r-\varrho)\big(D-\log(r-\varrho)\big)^{1/(1+\alpha)}\Big)\Big]^{\alpha}} \Bigg]
\end{aligned}
\right.
\end{equation}

By using $\log(1+t)/t\rightarrow1$ as $t\rightarrow0$ we find that for any fixed $\bar{C}$ we have 
\[h^{-1}(r)(r-\varrho)^{\alpha}\Bigg[\frac{1}{l_{1}}\log\Big(1+l_{1}\bar{C}(r-\varrho)\big(D-\log(r-\varrho)\big)^{1/(1+\alpha)}\Big)\Bigg]^{-\alpha}\bar{C}^{-\alpha}\rightarrow1\]
as $r\rightarrow\varrho.$ 
Now consider
\begin{equation}
\left\{
\begin{aligned}
&\frac{c_{2}(r-\varrho)^{\alpha}\Big(1+l_{1}(r-\varrho)\big(D-\log(r-\varrho)\big)^{1/(1+\alpha)}\Big)}{\bar{C}h(r)\Big[\frac{1}{l_{1}}\log\Big(1+l_{1}\bar{C}(r-\varrho)\big(D-\log(r-\varrho)\big)^{1/(1+\alpha)}\Big)\Big]^{\alpha}}
=\frac{c_{2}\Bigg(l_{1}(r-\varrho)\bar{C}\big(D-\log(r-\varrho)\big)^{1/(1+\alpha)}\Bigg)^{\alpha}}{\bar{C}^{1+\alpha}\Big[\log\Big(1+l_{1}\bar{C}(r-\varrho)\big(D-\log(r-\varrho)\big)^{1/(1+\alpha)}\Big)\Big]^{\alpha}}+\\
&\frac{c_{2}l_{1}^{\alpha}(r-\varrho)^{\alpha}\Big(l_{1}\bar{C}(r-\varrho)\big(D-\log(r-\varrho)\big)^{1/(1+\alpha)}\Big)}{\bar{C}h(r)\Big[\log\Big(1+l_{1}\bar{C}(r-\varrho)\big(D-\log(r-\varrho)\big)^{1/(1+\alpha)}\Big)\Big]^{\alpha}}.
\end{aligned}
\right.
\end{equation}
Observe that the first term on the right hand side converge to $\bar{C}^{-\alpha-1}$ and second term converges to $0$ in view of the above observation and the fact $t^{a}\log(t)\rightarrow0$ as $t\rightarrow0$ for any $a>0.$
We may choose $\bar{C}$ larger so that $c_{2}\bar{C}^{-1-\alpha}<\lambda/10(1+\alpha)$ and $(1/l_{2})\log(1+l_{2}u_{4}(R))\geq \sup _{x \in \Omega} u(x)$. Then we choose $\theta$ small enough such that the contribution from second term and remaining part of the last term in the right hand side of \eqref{B6} is less than $\lambda/5(1+\alpha).$ \\

In order to obtain the lower bound for $u$ again we proceed for any $x\in\Omega_{\varrho/2}$ by finding $y=y(x)\in\partial\Omega$ and $z\in\Omega$ such that $B_{\varrho}(z)$ is the ball touching $\partial\Omega$ at $y\in\partial\Omega.$ Now we want to apply the comparison principle in the annular domain $B_{\varrho}(z)\setminus B_{\varrho/2}(z).$ Let us consider the following function 
\[u_{5}(r)=\bar{c}(\varrho-r)(D-\log(\varrho-r))^{1/1+\mu},\]
where $D$ and $\bar{c}$ will be chosen later. Let us show that for $\bar{c}$ sufficiently small such that 
\begin{equation}\label{B666}
F^{-}_{2}(u_{5})+\frac{c_{0}}{l_{2}}(1-l_{2}u_{5})\log(1-l_{2}u_{5})+C_{1}(\varrho-r)^{\alpha-1}\Big[(1-l_{2}u_{5}(r))\Big]\Big[\frac{1}{l_2}\log\Big(1-l_{2}u_{5}(r)\Big)\Big]^{-\alpha}\geq 0~~\text{for~all}~~r\in[\varrho/2,\varrho].
\end{equation}
Observe that $u_{5}$ is non-negative and can be made as small we want by choosing $\overline{c}$ small. We also note that the term involving $c_{0}$ is non-positive. Let us proceed by introducing 
\[h(r)=\Big(D-\log(\varrho-r)\Big)^{\frac{-\mu}{1+\mu}}.\]
With the above introduction of $h$ in hand we get 
\begin{equation*}\label{B66}
\begin{aligned}
&u_{5}^{\prime}(r)=\overline{c}h(r)\Bigg[\frac{1}{(1+\mu)}-\Big(D-\log(\varrho-r)\Big)\Bigg]\\
&h^{\prime}(r)=-\frac{\mu}{(1+\mu)}\frac{h(r)}{\Big[D-\log(\varrho-r)\Big](\varrho-r)}\\
&\frac{(\varrho-r)^{\alpha}}{u_{5}^{\mu}}=\frac{h(r)}{\overline{c}^{\mu}(\varrho-r)}~~~\text{as}~~1+\alpha=\mu\\
&u_{5}^{\prime\prime}(r)=-\frac{\overline{c}h(r)}{(1+\alpha)(T-r)\Big[D-\log(\varrho-r)\Big]^{2}}\Bigg[\frac{\alpha}{(1+\alpha)}+\Big\{D-\log(\varrho-r)\Big\}\Bigg].\\
\end{aligned}
\end{equation*}
Taking into account the sign of $u^{\prime}_{5}(r),u^{\prime\prime}_{5}(r)$ we get the following:
\[\mathcal{M}_{\lambda,\Lambda}^-(D^2u)=\frac{\overline{c}h(r)\Lambda}{(\varrho-r)(1+\mu)}A_{\lambda,\Lambda}(r),\]
where
\[A_{\lambda,\Lambda}(r)=\Bigg[\frac{(\varrho-r)(N-1)}{r}-\frac{(1+\mu)(\varrho-r)(N-1)}{r}\Big[D-\log(\varrho-r)\Big]+\frac{1}{(1+\mu)\Big[D-\log(\varrho-r)\Big]^{2}}-(2+\mu)\Bigg].\]
Thus we find that
\begin{equation}\label{B66}
\left\{
\begin{aligned}
&F^{-}_{2}(u_{5})+-\frac{c_{0}}{l_{2}}(1-l_{2}u_{5})\log(1-l_{2}u_{5})+\frac{C_{4}(\varrho-r)^{\alpha}\Big[(1-l_{2}u_{5}(r))\Big]}{\Big[-\frac{1}{l_2}\log\Big(1-l_{2}u_{5}(r)\Big)\Big]^{\mu}}=\\
&\frac{\overline{c}h(r)}{((\varrho-r))}\Bigg[\frac{\Lambda A_{\lambda,\Lambda}(r)}{(1+\mu)}+b(\varrho-r)\Big[\frac{1}{(1+\mu)}-\Big(D-\log(\varrho-r)\Big)\Big]-\frac{c_{0}(\varrho-r)}{l_{2}\overline{c}h(r)}\Big(1-l_{2}u_{5}(r)\Big)\log(1-l_{2}u_{5})\\
&+\frac{C(\varrho-r)^{1+\alpha}(1-l_{2}u_{5})}{\overline{c}h(r)\Big[-\frac{1}{l_{2}}\log(1-l_{2}u_{5})\Big]^{\mu}}\Bigg].
\end{aligned}
\right.
\end{equation}
Observe that $u_{5}(r)=\overline{c}(\varrho-r)h^{\frac{-1}{\mu}}(r)$ and $1+\alpha=\mu.$ Now, consider the following last two terms from the right hand side of the above expression
\begin{equation}\label{B6}
\left\{
\begin{aligned}
&-\frac{c_{0}(\varrho-r)}{l_{2}\overline{c}h(r)}\Big(1-l_{2}u_{5}(r)\Big)\log(1-l_{2}u_{4})+\frac{C(\varrho-r)^{1+\alpha}(1-l_{2}u_{4})}{\overline{c}h(r)\Big[-\frac{1}{l_{2}}\log(1-l_{2}u_{5})\Big]^{\mu}}\\
&=-\frac{c_{0}(\varrho-r)}{l_{2}\overline{c}h(r)}\Big(1-l_{2}u_{5}(r)\Big)\log(1-l_{2}u_{5})+\frac{C}{\overline{c}^{\mu}}\Bigg[\frac{-l_{2}u_{5}}{\log(1-l_{2}u_{5})}\Bigg]^{\mu}\Bigg[\frac{1}{\overline{c}}-\frac{(\varrho-r)}{h(r)}\Bigg].
\end{aligned}
\right.
\end{equation}
Note that
\[\frac{(\varrho-r)}{h(r)}=(\varrho-r)\Big[D-\log(\varrho-r)\Big]^{\frac{\mu}{1+\mu}}\]
is bounded for $\frac{\varrho}{2}\leq r\leq \varrho.$ As we know $\frac{\log(1+x)}{x}$ approaches to $1$ as $x$ tends to $0.$ Therefore, we can choose $\overline{c}$ sufficiently small such that 
\[\frac{1}{2}\leq\frac{-l_{2}u_{5}(r)}{\log(1-l_{2}u_{5})}\leq 1.\]
Other terms in the braces on the right-hand side of \eqref{B66} are bounded, so by the above calculation, it is clear that we can choose $\overline{c}$ sufficiently small such that the right-hand side of \eqref{B66} is non-negative. Thus $u_{5}$ is a subsolution of \eqref{B666}. We decrease $\overline{c},$ if necessary, so that $(-1/l_{2})\log(1-l_{2}u_{5}(\gamma/2))\leq m,$ where 
\[
m = \min\{u(x) \mid d(x) \geq \rho/2\}.
\]

Therefore, we find that $(-1/l_{2})\log(1-l_{2}u_{5}(r))$ is a subsolution and by virtue of the choice of $\bar{c}$ it is less than $u$ on the boundary of the annular region $B_{\varrho}(z)\setminus B_{\varrho/2}(z).$ Consequently, the comparison principle and Lemma \ref{calc} implies the lower estimate.\\
\textbf{Case~$1+\mu<\alpha.$}\\
Again as in the above case in order to show that in order to show
\begin{equation}\label{ca2}
u(x) \leq A_2 d(x)(D-\log d(x))^{(\mu+2) /(1+\alpha)}
\end{equation}
for some constant $A_{2}$ it is sufficient to show that 
\begin{equation}\label{ca3}
u(|x|)\leq \frac{1}{l_{2}}\log(1+l_{2}u_{6}(r)),
\end{equation}
 where $u_{6}(r)=\bar{C}(r-\varrho)^{(2+\mu) /(1+\mu)},$ for sufficiently large $\bar{C}.$  Then \eqref{ca2} follows from \eqref{ca3} and Lemma \ref{calc}. Again in view of Lemma \ref{trans}, it is sufficient to show that $u_{6}$ satisfies 
 \begin{equation*}
 F^{+}_{2}(u_{6})-\frac{c_{0}}{l_{1}}(1+l_{1}u_{6})\log(1+l_{1}u_{6})+c_{2}(r-\varrho)^{\alpha-1}\Big[(1+l_{1}u_{6}(r))\Big]\Big[\frac{1}{l_1}\log\Big(1+l_{1}u_{6}(r)\Big)\Big]^{-\alpha}\leq 0~~\text{for~all}~~r\in[\varrho,R].
 \end{equation*}
 As $u_{6}$ is non-negative so it is sufficient to show that  
\begin{equation}\label{ca4}
F^{+}_{2}(u_{6})+c_{2}(r-\varrho)^{\alpha-1}\Big[(1+l_{1}u_{6}(r))\Big]\Big[\frac{1}{l_1}\log\Big(1+l_{1}u_{6}(r)\Big)\Big]^{-\alpha}\leq 0~~\text{for~all}~~r\in[\varrho,R]
\end{equation}
for $\bar{C}$ sufficiently large satisfying $ \displaystyle{\sup_{\Omega} u(x)\leq \frac{1}{l_{1}} \log(1 + l_{1} u_{6}(R))}.$
Now let us proceed to show \eqref{ca4} by computing 
\begin{equation}\label{B6}
\left\{
\begin{aligned}
&u_3(r)=\bar{C}(r-\varrho)^{(2+\mu) /(1+\mu)}\\
&u_3^{\prime}(r)=\frac{\bar{C}(2+\mu)}{1+\alpha}(r-\varrho) h(r)\\
&u_3^{\prime\prime}(r)=\frac{\bar{C}(2+\mu)(1+\mu-\alpha)}{(1+\alpha)^2} h(r)\\
\end{aligned}
\right.
\end{equation}
where $h(r)=(r-\varrho)^{(2+\mu) /(1+\alpha)-2}.$ Now observe that 
\begin{equation}
\begin{aligned}
\mathcal{M}_{\lambda, \Lambda}^{+}\left(D^2u_6\right)=\frac{\bar{C}(2+\mu) h(r)}{1+\alpha}\left(\frac{\lambda(1+\mu-\alpha)}{1+\alpha}+\frac{\Lambda(N-1)(r-\varrho)}{r}\right).
\end{aligned}
\end{equation}
Observe that $(r-\theta)/r\leq \theta$ for all $r\in (\varrho, \varrho)$ so by choosing $\theta$ small if needed we get the following 
\begin{equation}
\mathcal{M}_{\lambda, \Lambda}^{+}\left(D^2u_{6}\right) \leq \bar{C} \lambda \frac{(2+\mu)(1+\mu-\alpha)}{2(1+\alpha)^2} h(r)
\end{equation}
for all $r \in[\varrho, R]$.
\begin{equation}\label{B6}
\left\{
\begin{aligned}
&F^{+}_{2}(u_{6})+\frac{c_{2}(r-\varrho)^{\mu}\Big[(1+l_{1}u_{6}(r))\Big]}{\Big[\frac{1}{l_1}\log\Big(1+l_{1}u_{6}(r)\Big)\Big]^{\alpha}}\leq\\
&\bar{C}h(r)\Bigg[\frac{\lambda(2+\mu)(1+\mu-\alpha)}{2(1+\alpha)^{2}}+b\frac{(2+\mu)(r-\varrho)}{1+\alpha}+\frac{c_{2}(r-\varrho)^{\mu}\Big(1+l_{1}\bar{C}(r-\varrho)^{(2+\mu)/(1+\alpha)}\Big)}{\bar{C}h(r)\Big[\frac{1}{l_{1}}\log\Big(1+l_{1}\bar{C}(r-\varrho)^{(2+\mu)/(1+\alpha)}\Big)\Big]^{\alpha}}\Bigg].
\end{aligned}
\right.
\end{equation}
Let us rewrite the last term and use the expression for $h$
\begin{equation}\label{B6}
\left\{
\begin{aligned}
&\frac{c_{2}(r-\varrho)^{\mu}\Big(1+l_{1}\bar{C}(r-\varrho)^{(2+\mu)/(1+\alpha)}\Big)}{\bar{C}h(r)\Big[\frac{1}{l_{1}}\log\Big(1+l_{1}\bar{C}(r-\varrho)^{(2+\mu)/(1+\alpha)}\Big)\Big]^{\alpha}}=\frac{c_{2}\Big[\bar{C}l_{1}(r-\varrho)^{(\mu+2)/(1+\alpha)}\Big]^{\alpha}}{\bar{C}^{1+\alpha}\Big[\log\Big(1+l_{1}\bar{C}(r-\varrho)^{(2+\mu)/(1+\alpha)}\Big)\Big]^{\alpha}}\\
&+\frac{c_{2}l_{1}(r-\varrho)^{(1+\mu)/(1+\alpha)}\Big[l_{1}\bar{C}(r-\varrho)^{(\mu+2)/(1+\alpha)}\Big]^{\alpha}}{\bar{C}^{\alpha}\Big[\log\Big(1+l_{1}\bar{C}(r-\varrho)^{(2+\mu)/(1+\alpha)}\Big)\Big]^{\alpha}}.
\end{aligned}
\right.
\end{equation}
We first choose $\bar{C}$ large enough as in the above the above step, then choose $\theta$ small enough such that $u_{6}$ is a super solution.\\
In order to show the lower estimate  as in the above case we consider 
\begin{equation}
u_{7}(r)=\overline{c}(\varrho-r)^{\frac{(2+\mu)}{(1+\alpha)}}.
\end{equation}
Now by proceeding as in the above case we find that 
\begin{equation}\label{B444}
\left\{
\begin{aligned}
&u^{\prime}_{7}(r)=-\overline{c}\Bigg(\frac{2+\alpha}{1+\mu}\Bigg)(\varrho-r)^{(2+\alpha)/(1+\mu)-1}\\
&=-\overline{c}\Bigg(\frac{2+\alpha}{1+\mu}\Bigg)(\varrho-r)h(r)\\
&u_{7}^{\prime\prime}(r)=\overline{c}\Bigg(\frac{2+\alpha}{1+\mu}\Bigg)\Bigg(\frac{2+\alpha}{1+\mu}-1\Bigg)h(r)\\
&\mathcal{M}_{\lambda, \Lambda}^{-}\left(D^2u_7\right)=\overline{c}h(r)\Bigg[ (N-1)\Bigg(\frac{2+\alpha}{1+\mu}\Bigg)\frac{(\Gamma-r)}{r}+\Bigg(\frac{2+\alpha}{1+\mu}\Bigg)\Bigg(\frac{2+\alpha}{1+\mu}-1\Bigg)\Bigg]\\
&=\overline{c}h(r)A(r)
\end{aligned}
\right.
\end{equation}
where 
\begin{equation*}
A(r)=\Bigg[ (N-1)\Bigg(\frac{2+\alpha}{1+\mu}\Bigg)\frac{(\Gamma-r)}{r}+\Bigg(\frac{2+\alpha}{1+\mu}\Bigg)\Bigg(\frac{2+\alpha}{1+\mu}-1\Bigg)\Bigg].
\end{equation*}
\begin{equation}\label{B444}
\left\{
\begin{aligned}
&F^{-}_{2}(u_{7})+-\frac{c_{0}}{l_{2}}(1-l_{2}u_{7})\log(1-l_{2}u_{7})+\frac{C_{5}(\varrho-r)^{\alpha}\Big[(1-l_{2}u_{7}(r))\Big]}{\Big[-\frac{1}{l_2}\log\Big(1-l_{2}u_{7}(r)\Big)\Big]^{\mu}}=\\
&\overline{c}h(r)\Bigg[A_{\lambda,\Lambda}(r)-b(\varrho-r)\Bigg(\frac{2+\alpha}{1+\mu}\Bigg)-\frac{c_{0}}{l_{2}\overline{c}h(r)}\Big(1-l_{2}u_{7}(r)\Big)\log(1-l_{2}u_{7})+\frac{C_{5}(\varrho-r)^{\alpha}\Big(1+l_{2}\overline{c}(\varrho-r)^{(2+\alpha)/(1+\mu)}\Big)}{\overline{c}h(r)\Big[-\frac{1}{l_{2}}\log\Big(1+l_{2}\bar{C}(\varrho-r)^{(2+\alpha)/(1+\mu)}\Big)\Big]^{\mu}}\Bigg].
\end{aligned}
\right.
\end{equation}
Let us consider the last two term in the right hand side of the above expression
\begin{equation}\label{B444}
\left\{
\begin{aligned}
&-\frac{c_{0}}{l_{2}\overline{c}h(r)}\Big(1-l_{2}u_{7}(r)\Big)\log(1-l_{2}u_{7})+\frac{C_{5}(\varrho-r)^{\alpha}\Big(1+l_{2}\overline{c}(\varrho-r)^{(2+\alpha)/(1+\mu)}\Big)}{\overline{c}h(r)\Big[-\frac{1}{l_{2}}\log\Big(1+l_{2}\bar{C}(\varrho-r)^{(2+\alpha)/(1+\mu)}\Big)\Big]^{\mu}}\\
&-\frac{c_{0}}{l_{2}\overline{c}h(r)}\Big(1-l_{2}u_{7}(r)\Big)\log(1-l_{2}u_{7})+\frac{C_{5}}{\overline{c}^{\mu}}\Bigg[\frac{\log(1-l_{2}u_{7})}{-l_{2}u_{7}}\Bigg]^{-\mu}\Bigg[\frac{1}{\overline{c}}-l_{2}(\varrho-r)^{(2+\alpha)/(1+\mu)}\Bigg].
\end{aligned}
\right.
\end{equation}
Now by choosing $\overline{c}$ sufficiently small we find that $u_{7}$ is a sub olution and rest of the proof follows on the same line.
\end{proof}
\section{Global $\Gamma$ Regularity}

This section is devoted to the proof of the regularity of solutions to \eqref{main}. Depending on the relation between $\alpha$ and $\mu$, we show that the solutions are globally H\"{o}lder continuous or have globally H\"{o}lder continuous gradients. In order to prove these results, we need the analogous result in the interior case. Observe that:
\begin{equation}
\tilde{F}(x, r, p, X) = F(x, r, p, X) - c_{0}r
\end{equation}

satisfies $(SC)^{\mu}$ in \cite{nornberg2019c1} with $w(r)=|r|$ and $d(x)=d+c_{0}.$ More precisely,
\begin{equation}\label{scmu}
(SC)^{\mu}\left\{
\begin{aligned}
&(i)~~\mathcal{M}_{\lambda, \Lambda}^{-}(M-N)-B(|p|+|q|)|p-q|-b|p-q|-d|r-s| \\
&\leqq F(M,p,r,x)-F(N,q,s,x)\\
&\leq\mathcal{M}_{\lambda, \Lambda}^{+}(M-N)+B(|p|+|q|)|p-q|+b|p-q|+d|r-s|,\\
&(ii)~~F(x,0,0,0)=0.
\end{aligned}
\right.
\end{equation}

Therefore, we can apply 2\cite{sirakovsolvability} and 1.1\cite{nornberg2019c1} to obtain the the interior H\"{o}lder regularity and interior $C^{1,\alpha}$ estimate for the solutions of \eqref{main}. Therefore, from now onwards we will write $F$ in place of $F-c_{0}.$\\
Next we observe that the class of equations considered here are invariant under diffeomorphisms(see \cite{nornberg2018methods}) and $\partial\Omega\in C^{2}$ we only need to prove regularity and estimates for some half ball, say $B^{+}_{\Gamma_{0}}$ for some positive $\Gamma_{0}.$ Indeed, for each fixed $x_{0}\in\partial\Omega$ there exists $B_{\Gamma_{0}}$ and $C^{2}$ diffeomorphism $y=\Phi(x)$ such that $\Phi(B_{\Gamma_{0}}\cap\partial\Omega)$ is a hyperplane portion of the boundary of $\Phi(B_{\Gamma_{0}}\cap\Omega),$ say $\{y_{N}=0\}.$ By virtue of the above discussion it is sufficient to consider the equation in $B^{+}_{\Gamma_{0}},$ that is,
\begin{equation}\label{main3}
\left\{
\begin{aligned}
F(x,u,Du,D^{2}u)&=g(x^{\prime},x_{N})~\text{in}~B^{+}_{\Gamma_{0}}\\
u&=0~~\text{on}~~T,
\end{aligned}
\right.
\end{equation}
where $T=B^{+}_{\Gamma_{0}}\cap\{x_{N}=0\}$ and $F$ satisfies $(SC)$ with $\lambda,\Lambda,b, B,d $ possibly depending on $\Phi.$         
\subsection{H\"{o}lder~$\Gamma$ regularity}
\begin{theorem}\label{mainmain}
Let $F$ as above.
\begin{enumerate}
\item{} Assume $1+\mu<\alpha$ and \eqref{est1} holds. Then the solution of \eqref{main} is in $C^{\frac{\mu+2}{1+\alpha}}(\overline{\Omega}).$
\item{} Assume $1+\mu=\alpha$ and \eqref{est1} holds. Then the solution $u$ of \eqref{main} is in $C^{\beta}(\overline{\Omega}),$ for all $\beta<1.$
\end{enumerate}
\end{theorem}
\begin{proof}
\textbf{Case~$1+\mu<\alpha:$}~From the above discussion and $(SC)$ we know that the solution of \eqref{main} satisfies \eqref{main3}, and in view of Theorem \eqref{boundary1}(3), we may assume 
\begin{equation}
|g(x^{\prime},x_{N})|\leq Ax^{(\mu-\alpha)/(1+\alpha)}_{N}~\text{in}~B^{+}.
\end{equation}
Let $\Gamma_{2}<\Gamma_{1}<\Gamma_{0}$ and $T_{2}=B^{+}_{\Gamma_{2}}\cap\{x_{N}=0\}.$ There is $\gamma_{0}>0$ such that $(x^{\prime},x_{N})\in B^{+}_{\Gamma_{1}}$ for all $(x^{\prime},0)\in T_{2}.$ Then fix $(x^{\prime}_{0},0)\in T_{2},~\alpha_{1}=(2+\mu)/(1+\alpha)$ and $x_{\gamma}=(x^{\prime},3\gamma),$ and define the following scaled function
\begin{equation}
v(z)=\frac{u(\gamma~z+x_{\gamma})}{\gamma^{\alpha_{1}}},~~\text{for}~~y\in B_{3},
\end{equation}
for $0<\gamma<\gamma_{0}.$ Note that this function satisfies 
\begin{equation}
F_{\gamma}[v]:=F(D^{2}w, \gamma Dw, \gamma y+x_{\gamma})=g_{\gamma}(y)~~\text{in}~~~B_{3},
\end{equation}
where $g_{\gamma}=g(z\gamma+x_{\gamma})\gamma^{2-\alpha_{1}}.$ Note that in view of the choice of $\alpha_{1},$ we can find a constant $C$ independent of $\gamma$ such that $|g_{\gamma}(y)|\leq C,$ for $y\in B_{3}.$  Also notice that $F_{\gamma}$ satisfies all the assumptions of Theorem 1.1\cite{nornberg2019c1} for all $\gamma.$ Therefore by the interior elliptic estimate we find that $v\in C^{1}(\overline{B_{2}}).$ Thus,

\[|w(z_{1})-w(z_{2})|\leq C|z_{1}-z_{2}|,~\text{for~all}~~z_{1},z_{2}\in B_{2},\]
for some constant $C$ independent of $\gamma.$ Scaling back to $u$ we get 
\begin{equation}\label{eq4}
|u(x_{1})-u(x_{2})|\leq C|x_{1}-x_{2}|^{\alpha_1}~\text{for~all}~~x_{1},x_{2}\in B_{\gamma}(x_{\gamma}).
\end{equation}
Note that the above estimate is uniform in $\gamma\in(0,\gamma_{0})$ and $x^{\prime}_{0}\in T_{2}.$ Next as in \cite{felmer2012existence}, for a given $(x_{0}^{\prime},0)\in T_{2}$ we analyze the one dimensional function $w(t)=u(x_{0}^{\prime},t).$ We claim that 
\begin{equation}\label{eq2}
|w(t_{1})-w(t_{2})|\leq C|t_{1}-t_{2}|^{\alpha_{1}},~~\text{for~all~}~s_{1},s_{2}\in [0,~4\gamma_{0}),
\end{equation}
where the constant $C$ can be chosen independently of $(x^{\prime}_{0},0)\in T_{2}.$ To prove the claim we define the sequence $\gamma_i=(1 / 2)^i \gamma_0$ for $i \in \mathbb{N}.$ We first assume that $s_1>s_2>0$ and then consider the case $s_2=0.$ Let $i \leq j$ be indices such that $s_1 \in[2 \gamma_{i+1}, 2 \gamma_i]$ and $s_2 \in[2 \gamma_{j+1}, 2 \gamma_j].$ If $i=j$ then $\left(x_{0}^{\prime}, s_1\right),\left(x_{0}^{\prime}, s_2\right) \in B_{\gamma_{i+1}}\left(x_{\gamma_{i+1}}\right)$ and the claim follows by \eqref{eq4}. If $j=i+1$ then
$$
\frac{\left|w\left(s_1\right)-w\left(s_2\right)\right|}{\left|s_1-s_2\right|^\gamma} \leq \frac{\left|w\left(s_1\right)-w\left(2 \gamma_{i+1}\right)\right|}{\left|s_1-s_2\right|^\gamma}+\frac{\left|w(2 \gamma_j)-w\left(s_2\right)\right|}{\left|s_1-s_2\right|^\gamma} \leq 2 C,
$$
since $s_1-s_2 \geq s_1-2 \gamma_{i+1}$ and $s_1-s_2 \geq 2 \gamma_j-s_2$. If $i<j+1$ and noting that $s_1-s_2 \geq(1 / 2)^{i+1} \gamma_0=(1 / 2)^{i+1-k}\left(2 \gamma_k-2 \gamma_{k+1}\right)$ we have
\[\begin{aligned}
\frac{\left|w\left(s_1\right)-w\left(s_2\right)\right|}{\left|s_1-s_2\right|^\gamma} \leq & \frac{\left|w\left(s_1\right)-w\left(2 \gamma_{i+1}\right)\right|}{\left|s_1-s_2\right|^\gamma}+\sum_{k=i+1}^{k-j-1}\left\{\frac{\left|w\left(2 \gamma_k\right)-w\left(2 \gamma_{k+1}\right)\right|}{\left|s_1-s_2\right|^\gamma}\right\} \\
& +\frac{\left|w(2 \gamma_j)-w(s_2)\right|}{\left|s_1-s_2\right|^\gamma}
\leq 2 C+C\sum_{k=i+1}^{k-j-1}\left(\frac{1}{2}\right)^{\gamma(k-i-1)}\leq kC,
\end{aligned}
\]
where $k$ is independent of $i,j.$ As the above series converges so the claim follows in the case $s_{2}>0.$ In the case $s_2=0$, we use the continuity of $u$ to obtain that
$$
w\left(s_1\right)-w(0)=w\left(s_1\right)-w\left(2 \gamma_{i+1}\right)+\sum_{k=i+1}^{k-\infty}\left\{w\left(2 \gamma_k\right)-w\left(2 \gamma_{k+1}\right)\right\},
$$
from where we proceed as before, completing the proof of the claim.
Then we prove H\"{o}lder continuity or $u$ in all $\overline{{B}_{\Gamma_2}^{+}}$ as follows. Given $x=\left(x^{\prime}, x_N\right), y=$ $\left(y^{\prime}, y_N\right) \in \bar{B}_{\Gamma_2}^{\prime}$ we consider two cases:
\begin{enumerate}
\item{} If we have $|x-y|<x_N / 3$ (or $|x-y|<y_N / 3$ ) we just apply \eqref{eq4} in a ball containing both $x$ and $y$.
\item{} Otherwise by using, \eqref{eq2} we have

\[\begin{aligned}
& \frac{\left|u(x)-u\left(x^{\prime}, 0\right)\right|}{|x-y|^\gamma} \leq 3^\gamma \frac{\left|u(x)-u\left(x^{\prime}, 0\right)\right|}{x_N^\gamma} \leq 3^\gamma C, \\
& \frac{\left|u(y)-u\left(y^{\prime}, 0\right)\right|}{|x-y|^\gamma} \leq 3^\gamma \frac{\left|u(y)-u\left(y^{\prime}, 0\right)\right|}{x_N^\gamma} \leq 3^\gamma C, \\
&
\end{aligned}\]
\end{enumerate}
In view of $u\left(x^{\prime},0\right)=u\left(y^{\prime}, 0\right)=0$, H\"{o}lder continuity of $u$ follows.\\
In case $\mu-\alpha=-1$, we just take any $\alpha_1<1$ and the same argument applies
Next, we are going to show that in the case $\mu-\alpha+1>0$ we can improve the global regularity of the solution, obtaining a Holder estimate for the gradient in $\bar{\Omega}$.
\end{proof}

\subsection{\texorpdfstring{$C^{1,\mu}(\overline{\Omega})$}{C^{1,mu}(Omega)}:}
We proceed by establishing the gradient estimate on the boundary by following Krylov's approach. Such estimate for equations with quadratic growth in gradient has already been established in \cite{braga2020krylov}. But the approach in \cite{braga2020krylov} requires that $g\in L^{N}$ which, in general, is not true in this case. Therefore, we follow the approach mentioned in \cite{felmer2012existence} but we do not have Harnack inequality for such equations. We transform the equation by applying Lemma 2.3\cite{sirakovsolvability} and then apply the weak Harnack inequality. In order to state the first result we need the following notation. We denote by 
\[B_{R,\delta}=\{~(x^{\prime},x_{N})~~|~~|x^{\prime}|<R,~0<x_{N}<\delta R\}\]
and 
\[B^{\star}_{R/2,\delta}=\{~(x^{\prime},x_{N})~~~|~~~|x^{\prime}|<R,~~~\delta~R/2<x_{N}<3\delta~R/2\}.\]
\begin{proposition}\label{propkrylov}
There exists a $\delta>0$ such that if $u$ be a non negative strong supersolution of 
\begin{equation}\label{strong11}
\left\{
\begin{aligned}
\mathcal{M}_{\lambda,\Lambda}^-(D^2u)-B|Du|^2-b|Du|&\leq g(x^{\prime},x_{N})~~\text{in}~~B_{2R,\delta}\\
u&=0~~~\text{on}~~\{x_{n}=0\}\cap\overline{B_{2R,\delta}}.
\end{aligned}
\right.
\end{equation}
Then there exist two positive constants $p(N,\lambda,\Lambda,b)$ and $C_{2}(N,\lambda,\Lambda,b)$ such that
\begin{equation}\label{RTU1}
\Bigg[\fint_{B_{R / 2,\delta}^*}\Bigg(\frac{u}{x_{N}}\Bigg)^{p}dx\Bigg]^{1/p}\leq C_{2}e^{l_{1}\|u\|_{L^{\infty}(\Omega)}}\Bigg(\inf_{\substack{B_{R / 2,\delta}}}\Bigg(\frac{u}{x_{N}}\Bigg)+\|g\|_{L^{N}(B_{R / 2,\delta}^*)}\Bigg).
\end{equation}
\end{proposition}
\begin{proof}
We first prove that there is a positive $\delta$ such that
\begin{equation}\label{c1a2}
\inf_{\substack{\left|x^{\prime}\right|<R\\{x_N=\delta R}}}\Big(\frac{u}{x_{N}}\Big)\leq2\inf_{\substack{B_{R/2,\delta}}}\Big(\frac{u}{x_{N}}\Big).
\end{equation}
Since $u$ is non negative so without loss of generality we can assume that the left hand side is positive otherwise \eqref{c1a2} is true. 
Consider
\[U=\frac{1-e^{-l_{1}u }}{l_{1}}.\]
This function satisfies the following equation
\begin{equation}\label{TU1}
\begin{aligned}
\mathcal{M}_{\lambda,\Lambda}^-(D^2U)-b|DU|\leq
&e^{-l_{1}u }\Big(a_{2}x_{N}^{\mu-\alpha}\Big)\\
&\leq a_{2}x_{N}^{\mu-\alpha}.\\
\end{aligned}
\end{equation}
We have used that $u$ is non negative.
Now we consider the following equation 
\[Z(x^{\prime},x_{N})=\frac{1-e^{-l_{1}W(x^{\prime},x_{N})}} {l_{1}},\]
where
\[W(x^{\prime},x_{N})=L\Bigg(1-\frac{|x^{\prime}|^{2}}{R^{2}}+\frac{x^{1-\nu}_{N}-(R\delta)^{1-\nu}}{(R\delta)^{\frac{1-\nu}{2}}}\Bigg)x_{N}\]
and $L=\inf_{\substack{\left|x^{\prime}\right|<R\\{x_N=\delta R}}}\frac{u}{x_{N}}.$ Now we consider the following 
\begin{equation}\label{c1a1}
\begin{aligned}
\mathcal{M}_{\lambda,\Lambda}^-(D^2Z)-b|DZ|=&e^{\{-l_{1}W(x^{\prime},x_{N})\}}\big[\mathcal{M}_{\lambda,\Lambda}^+(D^2W-l_{1}DW\otimes DW)-b|DW|\big]\\
&\geq e^{\{-l_{1}W(x^{\prime},x_{N})\}}\big[\mathcal{M}_{\lambda,\Lambda}^-(D^2W)-B|DW|^{2}-b|DW|\big]\\
\end{aligned}
\end{equation}
Now we compute the expression in the bracket for $W.$
We can compute the derivative of the above function 
\[DW(x^{\prime},x_{N})=\frac{L}{R^{2}}\Bigg[-2 x_1 x_N \ldots,-2 x_{N-1} x_N,~~R^{2}-\left|x^{\prime}\right|^2~+~\frac{R^{2}\Big\{(2-\nu) x_N^{1-\nu}-(R\delta)^{1-\nu}\Big\}}{(R\delta)^{(1-\nu)/2}}\Bigg].\]
We find that the following estimate holds
\begin{equation*}
\begin{aligned}
|DW(x)|&\leq\frac{4L}{R^{2}}\sqrt{1+(2-\nu)^{2}\Bigg(\frac{x^{2}_{N}}{R\delta}\Bigg)^{1-\nu}+2(R\delta)^{1-\nu}}\\
&\leq\frac{4L}{R^{2}}\sqrt{1+5(R\delta)^{1-\nu}},
\end{aligned}
\end{equation*}
for all $x\in B_{2R,\delta}.$
Now consider the Hessian of $W,$
\[\begin{aligned}
& D^2 W(x)=\frac{L}{R^{2}}\left[\begin{array}{ccccc}
-2 x_N & 0 & \cdots & 0 & -2 x_1 \\
0 & -2 x_N & \cdots & 0 & -2 x_2 \\
\vdots &~ & ~\ddots & ~\vdots&\vdots \\
0 & \cdots & 0 & -2 x_N&-2 x_{N-1} \\
-2 x_1 & \cdots &\cdots& -2 x_{N-1} & \Gamma x_N^{-v}
\end{array}\right] \text {, } \\
&
\end{aligned}\]
with $\Gamma=(2-\nu)(1-\nu)/\left((R\delta)^{(1-\nu)/ 2}\right)$. Note that $\Gamma$ is large when $\delta$ is small. Now we compute the following
\begin{equation}\label{enarm}
\begin{aligned}
&\mathcal{M}_{\lambda,\Lambda}^-(D^2W) - l_{1}\lambda |DW|^{2} - b|DW|  \\
&\geq \frac{L}{R^{2}}\left[-2\Lambda(N-1)\delta R + \frac{\big\{R^{2}\lambda(2-\nu)(1-\nu)\big\}x^{-\nu}_{N}}{(\delta R)^{(1-\nu)/2}} \right.  
\left. + 2(\lambda - \Lambda)2R - 4b\sqrt{1+5(R\delta)^{1-\nu}}-\frac{16l_{1}L\Lambda}{R^{2}}\left(1+5(R\delta)^{1-\nu}\right) \right]  \\
\end{aligned}
\end{equation}

For a fixed value of $\nu\in(\alpha-\mu,~1)$ and $R$ we choose $\delta$ sufficiently small such that for any $x_{N}\in(0,\delta R)$ following holds:
\begin{equation*}\label{enarm}
\begin{aligned}
2\Lambda(N-1)\delta R+2(\Lambda-\lambda)2R+4b\sqrt{1+5(R\delta)^{1-\nu}}+\frac{16l_{1}L\Lambda}{R^{2}}\Bigg\{1+5(R\delta)^{1-\nu}\Bigg\}\leq\frac{\big\{R^{2}\lambda(2-\nu)(1-\nu)x^{-\nu}_{N}\big\}}{4(\delta R)^{(1-\nu)/2}},
\end{aligned} 
\end{equation*}
holds for any $x_{N}\in(0,\delta R).$ Furthermore, if necessary we can make it further smaller so that the following inequality also holds
\[a_{2}\leq \frac{L\lambda(2-\nu)(1-\nu)}{2(R\delta)^{\frac{1-\nu}{2}}}e^{-l_{1}\|W\|_{L^{\infty}}}.\]
Here we would like to remark that $W$ is bounded independent of $\delta.$ In fact, we can see that $\|W\|_{L^{\infty}(B_{2R,\delta})}\leq LR$ for any $\delta\in(0,1).$  Now, by equation \eqref{enarm},\eqref{c1a1} and above choice of $\delta$  we get
\begin{equation}
\mathcal{M}_{\lambda,\Lambda}^-(D^2Z)-b|DZ|>2a_{2}x_{N}^{-\nu}.
\end{equation}
Moreover we have the following 
\begin{enumerate}
\item[] On the bottom $x_{N}=0$ of $B_{R,\delta}=\left\{x~|~|x^{\prime}|<R, \quad 0<x_N<\delta R\right\},$ we have $u=W=0$ so $U=W$~~\text{on~the~bottom} of $B_{R,\delta}=\left\{x~|~|x^{\prime}|<R, \quad 0<x_N<\delta R\right\}.$
\item[] On $|x'|=R$, we have $W(x', x_{N}) \leq 0 \leq u(x', x_{N})$ and therefore $Z(x', x_{N}) \leq U(x', x_{N})$ on $|x'|=R$.
\item[] On the above portion $x_{N}=\delta R$ of $B_{R,\delta}=\left\{x~|~|x^{\prime}|<R, \quad 0<x_N<\delta R\right\},$ we have $W(x^{\prime},\delta R)\leq u(x^{\prime},R\delta),$ therefore $Z(x^{\prime},R\delta)\leq U(x^{\prime},R\delta).$
\end{enumerate}
Therefore, by the ABP maximum principle we have 
\[Z\leq U~~\text{on}~~B_{R,\delta}.\]
Consequently,
\begin{equation*}\label{enarm}
\begin{aligned}
\frac{1-e^{-l_{1}W(x^{\prime},x_{N})}} {l_{1}}\leq&\frac{1-e^{-l_{1}u(x^{\prime},x_{N})}} {l_{1}}\\
W(x^{\prime},x_{N})\leq&u(x^{\prime},x_{N}),\\
\end{aligned}
\end{equation*}
on $B_{R,\delta}.$ So we have 
\[L\leq2\inf _{\substack{B_{R/2,\delta}}}\Big\{\frac{u}{x_{N}}\Big\},\]
provided $(\delta R)^{(1-\nu)/2}\leq 1/4.$\\
Now note that $U$ is a non negative super solution to \eqref{TU1} in $B_{2R,\delta}$. Therefore, by Theorem \ref{weakh}, we can find two positive constants $p=p(N,\lambda, \Lambda, b)$ and $C=C(N,\lambda, \Lambda, b)$ such that 
\begin{equation}\label{enarm1}
\begin{aligned}
&\Big[\fint_{B_{R / 2,\delta}^*}U^{p}dx\Big]^{1/p}\leq C\Big(\inf_{\substack{B_{R / 2,\delta}^*}}U+R\|g\|_{L^{N}(B_{R}(0))}\Big)\\
\end{aligned}
\end{equation}
By using the following two inequalities
\begin{equation*}
\label{existence}
\left\{
\begin{aligned}
x&\geq 1-e^{-x}~~\implies u\geq U\\
x&\leq e^{x}-1~~~\implies e^{l_{1}u }U=\frac{e^{l_{1}u}-1}{l_{1}}\geq u
\end{aligned}
\right.
\end{equation*}
\eqref{enarm1} can be rewritten as follows 
\begin{equation}\label{TU2}
\begin{aligned}
&\Big[\fint_{B_{R / 2,\delta}^*}u^{p}dx\Big]^{1/p}\leq Ce^{l_{1}M}\Big(\inf_{\substack{B_{R / 2,\delta}^*}}u+R\|g\|_{L^{N}(B_{R}(0))}\Big).\\
\end{aligned}
\end{equation}
In $B_{R / 2,\delta}^*$ we have 
\begin{equation}\label{TU3}
\frac{2}{3\delta}\frac{u}{R}\leq\frac{u}{x_{N}}\leq\frac{2}{\delta}\frac{u}{R}.
\end{equation}
Dividing by $R$ in \eqref{TU2} and taking \eqref{TU3} into account we find that 
\begin{equation}\label{TU41}
\begin{aligned}
&\Big[\fint_{B_{R / 2,\delta}^*}\big(\frac{u}{x_{N}}\big)^{p}dx\Big]^{1/p}\leq Ce^{l_{1}M}\Big(\inf_{\substack{B_{R / 2,\delta}}}\big(\frac{u}{x_{N}}\big)+\|g\|_{L^{N}(B_{R / 2,\delta}^*)}\Big).\\
\end{aligned}
\end{equation}
Finally, \eqref{c1a2} and \eqref{TU41} gives \eqref{RTU1}.
\end{proof}
Note that $M_{2}=\max|u|$ so while we will be applying it to the problems we will take the maximum of $u_{1}$ and $u_{2}.$
\begin{remark}\label{remark11}
\begin{enumerate}
\item{} Note that the proof of Proposition \ref{propkrylov} also works if we take $\delta$ further small.
\item{} For any fixed $\delta>0$ we define 
\[M(2R)=\sup_{B_{2R.\delta}}v~~~\text{and}~~~m(2R)=\inf_{B_{2R,\delta}}v,\]
where $v=u/x_{N}$ and $u$ is a solution of \eqref{main}.
\item{} In view of the Lipschitz estimate Theorem \ref{boundary1} we have
\[A_{1}\leq M(2R),m(2R)\leq A_{2}.\]
In the above inequality $A_1$ and $A_{2}$ also depend on $\|u\|_{L^{\infty}(\Omega)},$ where $u$ is a solution of \eqref{main}.
\item{} Given a solution $u$ of \eqref{main} if we define 
\[u_{1}=M(2R)x_{N}-u~~\text{and}~~u_{2}=u-m(2R)x_{n},\]
then both the functions are non negative in $B_{2R,\delta}$ and for any $0<\delta<1$ we have
\begin{equation}\label{AA}
\|u_{i}\|_{L^{\infty}(\Omega)}\leq A_{2}+\|u\|_{L^{\infty}(\Omega)}~~~\text{for}~~~i=1,2.
\end{equation}
\end{enumerate}
\end{remark}
\begin{theorem}\label{bmain}
Suppose that $u$~ is a solution of \eqref{main3} where $F$ is same as in Theorem \ref{mainmain} and 
\begin{equation}
|g(x^{\prime},x_{N})|\leq Ax_{N}^{\mu-\alpha}~~\text{in}~~B^{+}.
\end{equation}
Then there are $R_{1}\in (0, R_{0}),$ ~$\alpha,$ ~$C$ such that for any $0<R\leq R_{1}$
\begin{equation}\label{osc11}
\text{osc}_{\substack{{B^{+}_{R}}}}\Big(\frac{u}{x_{N}}\Big)\leq CR^{\alpha}.
\end{equation}
\end{theorem}
\begin{proof}
We observe that from $W^{2,p}$ estimate in \cite{nornberg2019c1} we find that $u$ is a strong solution.
Now let us consider $u_{1}=Mx_{N}-u$ and $u_{2}=u-mx_{N}$ where $M=M(2R)$ and $m=m(2R)$ as in Remark \ref{remark11} (ii). Note that both are nonnegative in  $B_{2R,\delta}$ and belong to $W^{2,N}_{\text{loc}}.$
Moreover, $u_{1}$ and $u_{2}$ solve respectively
\begin{equation}
\left\{
\begin{aligned}
&F_{1}(x,u_{1},Du_{1},D^{2}u_{1})=-g(x^{\prime},x_{N})+f_{1}(x)\\
&F_{2}(x,u_{2},Du_{2},D^{2}u_{2}=g(x^{\prime},x_{N})+f_{2}(x)
\end{aligned}
\right.
\end{equation}
in $B_{2R,\delta},$ where
\begin{equation*}
\begin{aligned}
&f_{1}(x)=F(x, Me_{N},Mx_{n},0)~~\text{and}~~~f_{2}=-F(x,mx_{N},Me_{N},0)\\
&F_{1}(x,r,p,X)=-F(x, Mx_{N}-r, Me_{N}-p,-X)+F(x, Mx_{N}, Me_{N},0)\\
&F_{2}(x,r,p,X)=F(x,r+mx_{n},p+me_{N},X)-F(x,mx_{N}, Me_{N},0),
\end{aligned} 
\end{equation*}
where $M=M(2R)$ and $m=m(2R).$ Observe that $F_{1}$ and $F_{2}$ satisfies $(SC)^{\mu}$ with same $B,d$ and $\tilde{b}=b+2A_{2}$ as in \eqref{scmu}. Here, we would like to remark again that $A_{2}$ depends on $\|u\|_{L^{\infty}(\Omega)}.$ Let us define $g_{1}=-g+f_{1}+du_{1}$ and $g_{2}=g+f_{2}+du_{2}.$ 
Since $F$ satisfies $(SC)^{\mu},$ $M,m$ satisfies the inequality in Remark \ref{remark11}(iii) and by Theorem \ref{boundary1}(i) we have
\begin{equation*}\label{B444}
\left\{
\begin{aligned}
&|g_{1}(x^{\prime},x_{N})|\leq |g(x^{\prime},x_{N})|+|f_{1}|+d|u|+dM(2R)x_{N}\\
&\leq (A_{2}+BA^{2}_{2}+bA_{2}+3dA_{2})x^{\mu-\alpha}_N\\
&\leq C_{1}(A_{2},B,b,d)x^{-\nu}_{N}\\
&\text{similarly}\\
&|g_{2}(x^{\prime},x_{N})|\leq (A_{2}+BA^{2}_{2}+bA_{2}+3dA_{2})x^{\mu-\alpha}_N\\
&\leq C_{2}(A_{2},B,b,d)x^{-\nu}_{N},
\end{aligned}
\right.
\end{equation*}
where $\nu\in (\alpha-\mu,1)$ and $\delta$ is sufficiently small. Then by Proposition \ref{propkrylov} there exists $C_{3}(\lambda,\Lambda, N, A_{2},B),C_{4}(\lambda,\Lambda, N, A_{2},B)$ and $p_{1},p_{2}$ such that 
\begin{equation}
\left\{
\begin{aligned}
\Bigg[\fint_{B_{R / 2,\delta}^*}\Bigg(\frac{u_{1}}{x_{N}}\Bigg)^{p_{1}}dx\Bigg]^{1/p_{1}}&\leq C_{1}e^{l_{1}\Big(A_{2}+\|u\|_{L^{\infty}(\Omega)}\Big)}\Bigg(\inf_{\substack{B_{R / 2,\delta}}}\Bigg(\frac{u_{1}}{x_{N}}\Bigg)+\|g_{1}\|_{L^{N}(B_{R / 2,\delta}^*)}\Bigg).\\
\Bigg[\fint_{B_{R / 2,\delta}^*}\Bigg(\frac{u_{2}}{x_{N}}\Bigg)^{p_{2}}dx\Bigg]^{1/p_{2}}&\leq C_{2}e^{l_{1}\Big(A_{2}+\|u\|_{L^{\infty}(\Omega)}\Big)}\Bigg(\inf_{\substack{B_{R / 2,\delta}}}\Bigg(\frac{u_{2}}{x_{N}}\Bigg)+\|g_{2}\|_{L^{N}(B_{R / 2,\delta}^*)}\Bigg).
\end{aligned}
\right.
\end{equation}
Now we take $p=\min\{p_{1},p_{2}\}$ and $M(\frac{R}{2})=\sup_{B_{\frac{R}{2},\delta}}\Big(\frac{u}{x_{N}}\Big)$ and $m(\frac{R}{2})=\inf_{B_{\frac{R}{2},\delta}}\Big(\frac{u}{x_{N}}\Big)$ we find that 

\begin{equation*}
\begin{aligned}
&M(2R)-m(2R)=\Bigg[\fint_{B_{R / 2,\delta}^*}\Bigg(M(2R)-m(2R)\Bigg)^{p}dx\Bigg]^{1/p}
=\Bigg[\fint_{B_{R / 2,\delta}^*}\Bigg(\frac{M(2R)x_{N}-u+u-m(2R)x_{N}}{x_{N}}\Bigg)^{p}dx\Bigg]^{1/p}\\
&\leq C(p) \Bigg[\Big[\fint_{B_{R / 2,\delta}^*}\Bigg(\frac{u_{1}}{x_{N}}\Bigg)^{p_{1}}dx\Big]^{1/p_{1}}+\Big[\fint_{B_{R / 2,\delta}^*}\Big(\frac{u_{2}}{x_{N}}\Big)^{p_{2}}dx\Big]^{1/p_{2}}\Bigg]\\
&\leq C(p)\Bigg[C_{1}e^{l_{1}\Big(A_{2}+\|u\|_{L^{\infty}(\Omega)}\Big)}\Bigg(\inf_{\substack{B_{R / 2,\delta}}}\Bigg(\frac{u_{1}}{x_{N}}\Bigg)+\|g_{1}\|_{L^{N}(B_{R / 2,\delta}^*)}\Bigg)+C_{2}e^{l_{1}\Big(A_{2}+\|u\|_{L^{\infty}(\Omega)}\Big)}\Bigg(\inf_{\substack{B_{R / 2,\delta}}}\Bigg(\frac{u_{2}}{x_{N}}\Bigg)+\|g_{2}\|_{L^{N}(B_{R / 2,\delta}^*)}\Bigg)\Bigg]\\
&\leq C(p)\max\{C_{1},C_{2},C_{3},C_{4}\}e^{l_{1}\Big(A_{2}+\|u\|_{L^{\infty}(\Omega)}\Big)}\Bigg[\Big(M(2R)-m(2R)\Big)-\Big( M(\frac{R}{2})-m(\frac{R}{2})\Big) +R^{1-\nu}\Bigg].
\end{aligned}
\end{equation*}
Therefore, we have 
\[
O(2R) \leq \gamma \left( O\left(\frac{R}{2}\right) + CR^{1-\nu} \right)
\]
where $\gamma=(C-1)/C<1$ and $C=C(p)\max\{C_{1},C_{2},C_{3},C_{4}\}e^{l_{1}\Big(A_{2}+\|u\|_{L^{\infty}(\Omega)}\Big)}$ and $O(R)=M(R)-m(R).$ Now following the standard iteration process as in Lemma 8.23\cite{GT} we find a constant $C(\lambda,\Lambda,b,B,A_{1},A_{2}N,\|u\|_{L^{\infty}(\Omega)})$ and $\tau\in(0,1)$ such that \eqref{osc11} holds. 
\end{proof}
\begin{theorem}\label{mainmain}
Under the assumptions of Theorem \ref{bmain}, there are numbers $R_{1}(0,R_{0})$ and $\beta\in (0,1)$ such that $u\in C^{1,\beta}(~\overline{B}_{R_{1}})$
\end{theorem}
\begin{proof}
Theorem \ref{bmain} implies that the function $u$ is differentiable on $T_{1}=B_{R_{1}}\cap\{x_{N}=0\}$ and its gradient is H\"{o}lder continuous there. Moreover, 
\begin{equation}\label{uni11}
\lim _{\left(x^{\prime}, x_N\right) \rightarrow\left(x^{\prime}, 0\right)} \frac{\partial u\left(x^{\prime}, 0\right)}{\partial x_N}-\frac{u\left(x^{\prime}, x_N\right)}{x_N}=0
\end{equation}
uniformly in $(x^{\prime}_{0},0)\in T_{2},$ for any fixed $0<R_{2}<R_{1}.$ Now as in the proof of the H\"{o}lder continuous case, choose $\gamma_{0}$ such that $(x^{\prime},6\gamma_{0})\in B^{+}_{R_{{1}}}$ for any $(x^{\prime},0)\in T_{2}.$ At this point notice that $A(x^{\prime}):=\frac{\partial u\left(\bar{x}^{\prime}, 0\right)}{\partial x_N}$ is uniformly bounded by $\max\{A_{1},A_{2}\}.$ Now, for a fixed $(x^{\prime}_{0},0)\in T_{2},$ $\alpha_{1}=\min\{1-\nu,\alpha,3/4\}$ we define the following scaled function 
\[v(z)=\frac{1}{\gamma^{1+\alpha_{1}}}\big\{u(x_{\gamma}+\gamma z)-A(x^{\prime}_{0})\gamma(z_{N}+3)\big\} ~~\text{for}~~~y\in B_{3},\]
where $\gamma\in(0,\gamma_{0})$ and $x_{\gamma}=(x^{\prime}_{0},3\gamma).$ The function $v$ is bounded in $B_{3}$ because 
\[u(x)-x_N A_0(x^{\prime}_{0})\left(\bar{x}^{\prime}\right) \leq C x_N^{1+\alpha},\]
by Theorem \ref{bmain}. In addition to this this $v$ satisfies 
\[F_{\gamma}(y,v(y),Dv(y),D^{2}v(y))=\gamma^{1-\alpha_{1}}\tilde{g}(x_{\gamma}+\gamma y)~\text{in}~B_{3},\]
where \[F_{\gamma}(y,r,p,X)=\gamma^{1-\alpha_{1}}\Bigg[F(x_{\gamma}+\gamma~y,\gamma^{1+\alpha_{1}}r+A_{0}\gamma(y_{N}+3), \gamma^{\alpha_{1}}p+A_{0}e_{N},\gamma^{\alpha_{1}-1}X)-F(x_{\gamma}+3y,A_{0}\gamma(y_{N}+3),A_{0}e_{N},0)\Bigg].\]
and \[\tilde{g}(x_{\gamma}+\gamma y)=g(x_{\gamma}+\gamma y)-F(x_{\gamma}+\gamma y,A_{0}\gamma(y_{N}+3),A_{0}e_{N},0)\]
Above operator satisfies the Structure condition $(SC)^{\mu}$ given by \eqref{scmu} with $B, b$ and $d$ replaced by $\gamma^{1-\alpha_{1}}B, \gamma[b+2A_{0}B]$ and $\gamma^{2}d.$ At this point we would like to point out that $|A_{0}|\leq \max\{A_{1},A_{2}\},$ where $A_{1}$ and $A_{2}$ is from Theorem \ref{boundary1}(i). Thus, by virtue of the choice $0<\alpha_{1}<1$ and uniform bound on $A_{0}$ we find that for any $0<\gamma<1$ and $(\overline{x},0)\in T_{2},$ $F_{\gamma}$ satisfies the $(SC)^{\mu}$ with $B,$$b$ and $d$ replaced by $B,$ ~$b+2B\max\{A_{1},A_{2}\}$ and $d$ which are independent of $\gamma$ and $(\overline{x},0)\in T_{2}.$ Now we consider the term
\begin{equation*}
\begin{aligned}
&\Big|\gamma^{1-\alpha_{1}}\tilde{g}(x_{\gamma}+\gamma y)\Big|\leq \gamma^{1-\alpha_{1}} \Big|g(x_{\gamma}+\gamma y)\Big|+\gamma^{1-\alpha_{1}}\Big|F\Big(x_{\gamma}+\gamma y, A_{0}\gamma(y_{N}+3), A_{0}e_{N},0\Big)\Big|\\
&\leq \tilde{C}+\gamma^{1+\alpha_{1}}\Big[6d\gamma+b+B\max\{A_{1},A_{2}\}\Big]\max\{A_{1},A_{2}\}\leq C.
\end{aligned}
\end{equation*}
As $\|v\|_{L^{\infty}(B_{3})}\leq C$ so we can apply Theorem 1.1\cite{nornberg2019c1}, to get $\beta\in (0,1)$ and a constant $C$ such that 
\begin{equation}
|Dv(z_{1})-Dv(z_{2})|\leq C|z_{1}-z_{2}|^{\beta}~~~\text{for~all}~~z_{1},z_{2}\in B_{2}.
\end{equation}
Here $\beta$ and $C$ may depend on $\|u\|_{L^{\infty}}$ through $A_{1},A_{2}$ but independent on $(\overline{x},0)\in T_{2}$ and $\gamma.$ We decrease $\beta$,if necessary, so that $\beta\leq \alpha_{1}$ and we see that, by the definition of $v,$ we have 
\begin{equation}\label{hol}
|Du(x_{1})-Du(x_{2})|\leq C|x_{1}-x_{2}|^{\tau}~~\text{for~all}~x_{1},x_{2}\in B_{2}(x_{\gamma}),
\end{equation}
 where the estimate is uniform in $\gamma\in(0, \gamma_{0})$ and $(\overline{x},0)\in T_{2}.$ Now as in the proof of Theorem \ref{mainmain} or in the proof of Proposition 3\cite{felmer2012existence} we can find a constant $C$ such that for any $(\overline{x},0)\in T_{2}$ the one dimensional function $z(\cdot)=Du(\overline{x},\cdot)$ satisfies 
 \begin{equation}\label{uniform}
 \big|z(s_{1})-z(s_{2})\big|\leq C|s_{1}-s_{2}|^{\tau}~~\text{for}~~s_{1},s_{2}\in ( 0, 4\gamma_{0}].
 \end{equation}
The above result can be extended to hold $[0,4\gamma_{0}]$ as in the proof of Theorem \ref{mainmain} provided gradient of $u$ is continuous in  $B^{+}_{R_{2}}.$ Which is indeed the case. In fact, in view of \eqref{uniform} and mean value inequality for $u(x^{\prime},s)$ we find an $\chi\in (0,x_{N})$ such that 
\begin{equation}\label{uni13}
\begin{aligned}
& \Big|Du(x^{\prime},x_{N})-\frac{u(x^{\prime},x_{N})}{x_{N}}\Big|=|Du(x^{\prime},x_{N})-Du(x^{\prime},\chi)|^{\tau}\\
&\leq C|x_{N}-\chi|^{\tau}\leq C|x_{N}|^{\tau}.
\end{aligned}
\end{equation}
Then, by rewriting 
\begin{equation*}\label{uni14}
Du(x^{\prime},x_{N})-Du(\overline{x}^{\prime},0)=Du(x^{\prime},x_{N})-\frac{u(x^{\prime},x_{N})}{x_{N}}+\frac{u(x^{\prime},x_{N})}{x_{N}}-Du(\overline{x}^{\prime},0)
\end{equation*}
and taking \eqref{uni11} and \eqref{uni13} into account we find that gradient of $u$ is continuous. Finally, we prove holder continuity of the gradient on $\overline{B}_{R_{2}}$ as follows in standard way. In fact, for any $x,y\in \overline{B}_{R_{2}}$ we have the following two cases:
{\color{red} Change}
\begin{enumerate}
\item[] If we have $|x-y| < x_{N}/3$ or $|x-y| < y_{N}/3$ in this case H\"{o}lder continuity of the gradient follows from \eqref{hol}.  
\item[] Otherwise by \eqref{uniform} we have
\begin{equation*}
\begin{aligned}
&\frac{|Du(x^{\prime},x_{N})-Du(x^{\prime},0)|}{|x-y|^{\tau}}\leq 3^{\tau}\frac{|Du(x)-Du(x^{\prime},0)|}{x^{\tau}_{N}}\leq 3^{\tau}C\\
&\frac{|Du(y^{\prime},x_{N})-Du(y^{\prime},0)|}{|x-y|^{\tau}}\leq 3^{\tau}\frac{|Du(y)-Du(y^{\prime},0)|}{x^{\tau}_{N}}\leq 3^{\tau}C\\
&\text{and}\\
&\frac{|Du(x^{\prime},0)-Du(y^{\prime},0)|}{|x-y|^{\tau}}\leq 3^{\tau}\frac{|Du(x^{\prime},0)-Du(y^{\prime},0)|}{|x^{\prime}-y^{\prime}|^{\tau}}\leq C,\\
\end{aligned}
\end{equation*}
where we have used the H\"{o}lder continuity of $Du(\cdot,0)$ on $T_{2}$(which is a consequence of Proposition \ref{propkrylov}) and $|x-y|\geq |x^{\prime}-y^{\prime}|.$

\end{enumerate}

\end{proof}

\bibliography{RegWithSingNon.bib}
\bibliographystyle{abbrv}
\end{document}